\documentclass[12pt,letterpaper]{article}
\usepackage{amsbsy}
\usepackage{amscd}
\usepackage{amssymb,amsmath}
\usepackage{amsfonts}
\usepackage[english,russian]{babel}
\usepackage{authblk}
\usepackage{anysize}
\marginsize{2.5cm}{2cm}{1cm}{1cm}

\title{\bf{Applications of algebraic methods in solving the center-focus
problem}}
\author{Mihail\,Popa$^1$, \ Victor\,Pricop$^{1,2}$}
\affil{$^1$ Institute of Mathematics and Computer Sciences,
\\ Academy of Sciences of Moldova}
\affil{ $^{1,2}$ "Ion Creang\u{a}" State Pedagogical University,
Chi\c{s}in\u{a}u}
\date{\small{\textbf{E-mail}:
mihailpomd@gmail.com, pricopv@mail.ru}}
\begin{document}
\maketitle
\selectlanguage{english}
 \centerline{{\bf Abstract}}
\bigskip

The nonlinear differential system $
\dot{x}=\sum_{i=0}^{\ell}P_{m_i}(x,y),\
\dot{y}=\sum_{i=0}^{\ell}Q_{m_i}(x,y)$ is considered, where
$P_{m_i}$ and $Q_{m_i}$ are homogeneous polynomials of degree
$m_i\geq 1$ in $x$ and $y$, $m_0=1$. The set
$\{1,m_i\}_{i=1}^{\ell}$ consists of a finite number $(\ell<\infty)$
of distinct natural numbers. It is shown that the maximal number of
algebraically independent focal quantities that take part in solving
the center-focus problem for the given differential system with
$m_0=1$, having at the origin of coordinates a singular point of the
second type (center or focus), does not exceed
$\varrho=2(\sum_{i=1}^{\ell}m_i+\ell)+3.$  We make an assumption
that the number $\omega$ of essential conditions for center  which
solve the center-focus problem for this  differential system  does
not exceed  $\varrho$, i.\,e. $\omega\leq\varrho$.
\\


\section{Introduction}

\ \ \ \  The nonlinear differential system
\begin{equation}
\begin{gathered}\label{eq1}
\frac{dx}{dt}=\sum_{i=0}^{\ell}P_{m_i}(x,y),\
\frac{dy}{dt}=\sum_{i=0}^{\ell}Q_{m_i}(x,y)
\end{gathered}
\end{equation}
is considered, where $P_{m_i}$ and $Q_{m_i}$ are homogeneous
polynomials of degree $m_i\geq 1$ in $x$ and $y$, $m_0=1$. The set
$\{1,m_1,m_2,...,m_{\ell}\}$ consists of a finite number
$(\ell<\infty)$ of distinct natural numbers. The coefficients and
variables in polynomials $P_{m_i}$ and $Q_{m_i}$ take values from
the field of the real numbers $\mathbb{R}$.

It is known that if the roots of characteristic equation of the
singular point $O(0,0)$ of the system \eqref{eq1} are imaginary,
then the singular point $O$ is a center (surrounded by closed
trajectories) or a focus (surrounded by spirals) [1,5]. In this case
the origin of coordinates is a singular point of the second type.

Hereafter we denote the system \eqref{eq1} by
$s(1,m_1,m_2,...,m_{\ell})$.

The center-focus problem can be formulated as follows: {\it Let for
the system} $s(1,m_1,$ $m_2,...,m_{\ell})$ {\it the origin of
coordinates be a singular point of the second type} ({\it center or
focus}). {\it Find the conditions which distinguish center from
focus}. This problem was posed by H.\,Poincar\'{e} [1,2]. The basic
results were obtained by  A.\,M.\, Lyapunov [5]. It was shown that
the conditions for center are the vanishing of an infinite sequence
of polynomials (focal quantities)
\begin{equation}
\begin{gathered}\label{eq2}
L_1,L_2,...,L_k,...
\end{gathered}
\end{equation}
in the coefficients of right side of the system \eqref{eq1}. If at
least one of the quantities \eqref{eq2} is not zero, then the origin
of coordinates for the system \eqref{eq1} is a focus. These
conditions are necessary and sufficient.

In the case of the system \eqref{eq1} from Hilbert's theorem on the
finiteness of basis of polynomial ideals it follows that in the
mentioned sequence \eqref{eq2} only a finite number of conditions
for center {\it are essential}, the rest are consequences of them.
Then the center-focus problem for the system \eqref{eq1} takes the
following formulation: {\it How many polynomials} ({\it essential
conditions for center})
\begin{equation}
\begin{gathered}\label{eq3}
L_{n_1},L_{n_2},...,L_{n_{\omega}},...\ (n_i\in\{1,2,...,k,...\};\
i=\overline{1,\omega};\ \omega<\infty)
\end{gathered}
\end{equation}
{\it from \eqref{eq2} must be equal to zero in order that all other
polynomials \eqref{eq2} would vanish}?

The problem of determining the number $\omega$ of essential
conditions for center \eqref{eq3} is complicated. It is completely
solved for the systems $s(1,2)$ and $s(1,3)$ [8,11], for which we
have respectively $\omega=3$ and $5$. Until now $\omega$ has not
been known for the system $s(1,2,3)$. There exists only a
Zol\^{a}dek hypothesis, which is based mostly on intuition, that for
the system $s(1,2,3)$ the number $\omega\leq 13$. To the present day
this hypothesis has not been disproved. But in [12] it is proved
that for the system $s(1,2,3)$ $12$ focal quantities are not enough
for solving the center-focus problem in the complex plane.

It is natural to ask why there is still no answer about the value of
$\omega$ from \eqref{eq3} for any system
$s(1,m_1,m_2,...,m_{\ell})$?

We can explain this failure as follows: searching for a finite
$\omega$ from \eqref{eq3}, till now the researchers have used
basically a known approach, i.\,e. with the help of certain
calculations they constructed the explicit form of the first focal
quantities from \eqref{eq2}, without knowing a priori the number
$\omega$. Sometimes the existence of some geometric properties for
the system \eqref{eq1} was assumed, for example, the existence of
integral straight lines, conics and other curves. Then with their
help the attempts were made to show that the vanishing of the
available quantities implies the vanishing of other members of the
sequence \eqref{eq2}, often there was only a vague idea about their
expressions.

This approach gave quite unsatisfactory results. One of the reasons
is due to the enormous computing for focal quantities, which can not
be overcome using supercomputers even for the system $s(1,2,3)$, not
to mention more complicated systems $s(1,m_1,m_2,...,m_{\ell})$.
Therefore, it is clear that the results obtained in this direction
refer more to the systems \eqref{eq1} of special forms.

From what has been said above the following conclusion can be drawn:
 solving the center-focus problem is equivalent to finding the
essential conditions for center \eqref{eq3}, that requires knowledge
of the number $\omega$, the finiteness of which follows from
Hilbert's theorem on the finiteness of basis of polynomial ideals.

Therefore, the problem of finding the number $\omega<\infty$ {\it or
obtaining for it an argued numerical upper bound} ({\it even as a
hypothesis}), {\it which is still absent}, is a very important
condition of the complete solving of the center-focus problem for
the system \eqref{eq1}.

The last affirmation can be considered as a generalized center-focus
problem for the systems $s(1,m_1,m_2,...,m_{\ell})$, and obtaining
an answer to it will be qualified as perhaps one of the sufficient
conditions in solving the mentioned problem.

\section{Graded algebras of comitants of the system (1)}

\ \ \ \ In [5,6,7] the type of center-affine polynomial comitant
with respect to the center-affine group $GL(2,\mathbb{R})$ for any
differential system $s(m_0,m_1,m_2,...,m_{\ell})$ was determined,
and it is denoted as follows:
\begin{equation}
\begin{gathered}\label{eq4}
(d)=(\delta,d_0,d_1,...,d_{\ell}),
\end{gathered}
\end{equation}
where \ $\delta$ \ is \  the\  degree\  of \  homogeneity of this
comitant in phase \ variables \ $x,y$, and\  $d_i$\
$(i=\overline{1,\ell})$ is the degree of homogeneity of the same
comitant in the coefficients of the polynomials $P_{m_i}(x,y),\
Q_{m_i}(x,y)$ from the right side of the system \eqref{eq1}.

In [7] the following affirmations were proved:

\noindent\textbf{Proposition 1}. {\it The set of center-affine
comitants of the system \eqref{eq1} of the same type \eqref{eq4}
forms a finite linear space $V_{m_0,m_1,m_2,...,m_{\ell}}^{(d)}$,
i.\,e. it has a finite maximal system of linearly independent
comitants of the given type} ({\it linear basis}), {\it all the rest
are linearly expressed trough them}.

\noindent\textbf{Proposition 2}. {\it In order that any homogeneous
polynomial of the type \eqref{eq4} in phase variables and
coefficients of the system \eqref{eq1} would be a center-affine
comitant of this system, it is necessary and sufficient that it be
an unimodular comitant} ({\it invariant polynomial with respect to
the unimodular group $SL(2,\mathbb{R})$}) {\it of the same type
\eqref{eq4} for the given system}.

\noindent\textbf{Proposition 3}. [6] {\it For any center-affine
comitant of differential system of the type \eqref{eq4} the
following equality holds:
\begin{equation}
\begin{gathered}\label{eq501}
2g=\sum_{i=0}^{\ell}d_i(m_i-1)- \delta,
\end{gathered}
\end{equation}
where $g$ is usually called the weight of given comitant, and it is
an integer number}.

Following Propositions 1--2 and according to [7] we denote the space
of unimodular comitants of the type \eqref{eq4} for the system
$s(1,m_1,m_2,...,m_{\ell})$  by
\begin{equation*}
\begin{gathered}
S_{1,m_1,m_2,...,m_{\ell}}^{(d)}\cong
V_{1,m_1,m_2,...,m_{\ell}}^{(d)}.
\end{gathered}
\end{equation*}

Let us consider the linear space
\begin{equation}
\begin{gathered}\label{eq5}
S_{1,m_1,m_2,...,m_{\ell}}=\sum_{(d)}S_{1,m_1,m_2,...,m_{\ell}}^{(d)},
\end{gathered}
\end{equation}
which is a graded algebra of comitants of the system
$s(1,m_1,m_2,...,m_{\ell})$, where its components satisfy the
inclusion
\begin{equation*}
\begin{gathered}
S_{1,m_1,m_2,...,m_{\ell}}^{(d)}S_{1,m_1,m_2,...,m_{\ell}}^{(e)}\subseteq
S_{1,m_1,m_2,...,m_{\ell}}^{(d+e)},\
S_{1,m_1,m_2,...,m_{\ell}}^{(0)}=\mathbb{R}.
\end{gathered}
\end{equation*}

We denote by $SI_{1,m_1,m_2,...,m_{\ell}}$  a graded algebra of
unimodular invariants (comitants that do not depend on the phase
variables $x,y$) of the system $s(1,m_1,m_2,...,m_{\ell})$, which
satisfies the inclusion
\begin{equation}
\begin{gathered}\label{eq6}
SI_{1,m_1,m_2,...,m_{\ell}}\subset S_{1,m_1,m_2,...,m_{\ell}}.
\end{gathered}
\end{equation}

As for the first time the comitants and invariants for systems of
the form \eqref{eq1} were introduced by K.\,S.\, Sibirsky [14],
hereafter we will refer to these and similar algebras as {\it
Sibirsky algebras}.

\section{Krull dimension for Sibirsky graded algebras}

\ \ \ \ From the theory of invariants and tensors [5,13] it results
that the Sibirsky graded algebras $S_{1,m_1,m_2,...,m_{\ell}}$ and
$SI_{1,m_1,m_2,...,m_{\ell}}$ are commutative and finitely
determined algebras. If for these algebras we introduce a unified
notation $A$, then the last affirmation can be written as
\begin{equation}
\begin{gathered}\label{eq7}
A=<a_1,a_2,...,a_m\mid f_1=0,f_2=0,...,f_n=0>\ (m,n<\infty),
\end{gathered}
\end{equation}
where $a_i$ are generators for this algebra, and $f_j$ are defining
relations (syzygies).

It is known from [7] that for the simplest differential system
$s(0,1)$ of the form
\begin{equation}
\begin{gathered}\label{eq8}
\dot{x}= a+cx+dy,\ \dot{y}=b+ex+fy
\end{gathered}
\end{equation}
the finitely defined graded algebras of comitants $S_{0,1}$ and
invariants $SI_{0,1}$ can be written respectively
\begin{equation*}
\begin{gathered}
S_{0,1}=<i_1,i_2,i_3,k_1,k_2,k_3\mid(i_1k_1-k_3)^2+k_3^2-i_2k_1^2-2i_3k_2=0>,
\end{gathered}
\end{equation*}
\begin{equation}
\begin{gathered}\label{eq9}
SI_{0,1}=<i_1,i_2,i_3>,
\end{gathered}
\end{equation}
where
\begin{equation}
\begin{gathered}\label{eq10}
i_1=c+f,\ i_2=c^2+2de+f^2,\
i_3=-ea^2+(c-f)ab+db^2,\\ k_1=-bx+ay,\ k_2=-ex^2+(c-f)xy+dy^2,\\
k_3=-(ea+fb)x+(ca+db)y.
\end{gathered}
\end{equation}
We note that using the system \eqref{eq8} the whole theory of
center-affine (unimodular) comitants and invariants for
two-dimensional polynomial differential systems can be illustrated.

\noindent\textbf{Definition 1}. [15] {\it Elements $a_1,a_2,...,a_r$
of the algebra $A$ are called algebraically independent if for any
non-trivial polynomial $F$ in these $r$ elements the following
inequality holds}:
\begin{equation*}
\begin{gathered}
F(a_1,a_2,...,a_r)\neq0.
\end{gathered}
\end{equation*}

\noindent\textbf{Definition 2}. {\it The maximal number of
algebraically independent elements of an graded algebra $A$ is
called the Krull dimension of this algebra and is denoted by
$\varrho(A)$}.

It is known [7] that for an algebra $A$ of the form \eqref{eq7} the
equality $n=m-\varrho(A)$ holds. However, this equality is not very
effective because it is impossible to determine the numbers $m$ and
$n$ for most algebras of invariants and comitants for systems of the
form \eqref{eq1}.

In the classical theory of invariants [16] a set of elements
$a_1,a_2,$ $...,a_{\varrho(A)}$ from $A$ which define the Krull
dimension of the algebra $A$ is called {\it an algebraic basis}.
This means that for any $a\in A\ (a\neq a_j)$ there exists a natural
number $p$ such that the following identity holds:
\begin{equation}
\begin{gathered}\label{eq11}
P_0a^p+P_1a^{p-1}+...+P_p=0,
\end{gathered}
\end{equation}
where $P_k\ (k=\overline{0,p})$ are polynomials in $a_j\
(j=\overline{1,\varrho(A)})$. We note that in general
$P_0\not\equiv1$.

If for any $a\in A$ in \eqref{eq11} we have $P_0\equiv1$, then this
basis is called {\it integer algebraic basis}. The existence a basis
was shown by D.\,Hilbert (see [16]). We denote the number of its
elements by $\varrho^{\prime}(A)$.

We note that in general the numbers of elements in the mentioned
bases does not always coincide. For example, from [7] we have that
for the system $s(4)$ the Krull dimension $\varrho(SI_4)=7$, but
from [17] for the same system we obtain that the number of elements
in the integer algebraic basis of the same algebra is
$\varrho^{\prime}(SI_4)=9$, i.\,e.
$\varrho(SI_4)<\varrho^{\prime}(SI_4)$. From [7] we have that for
the system $s(0,1)$ the equality
$\varrho(S_{0,1})=\varrho^{\prime}(S_{0,1})=5$ holds, and
$\varrho(SI_{0,1})=\varrho^{\prime}(SI_{0,1})=3$. Also from [5,6,7]
and [18] it follows that for the systems $s(2)$ and $s(3)$ we have
$\varrho(SI_{2})=\varrho^{\prime}(SI_{2})=3$,
$\varrho(SI_{3})=\varrho^{\prime}(SI_{3})=5$. From [7] and [19] we
find that for the system $s(1,2)$ the equalities
$\varrho(SI_{1,2})=\varrho^{\prime}(SI_{1,2})=7$ are valid.
However, for the system $s(1,2,3)$ according to [7,20] we have that
$\varrho(SI_{1,2,3})=15$, but $\varrho^{\prime}(SI_{1,2,3})=21$.

The mentioned examples lead us to the relation
\begin{equation*}
\begin{gathered}
\varrho(A)\leq\varrho^{\prime}(A).
\end{gathered}
\end{equation*}

This inequality accentuates that the integer algebraic basis
contains an algebraic basis of an algebra $A$. The proof of this
fact can be easily obtained by an indirect proof.

\noindent\textit{Remark} 1.  The main property of an integer
algebraic basis of an algebra $A$ of invariants is that it is the
minimum number of elements of the algebra $A$ such that if they are
equal to zero, all elements of the algebra $A$ vanish.

Hereafter we need some evident affirmations:

\noindent\textbf{Proposition 4}. {\it If $B$ is a graded subalgebra
of an algebra $A$, then between the Krull dimensions of these
algebras the following inequality holds}:
\begin{equation*}
\begin{gathered}
\varrho(B)\leq\varrho(A).
\end{gathered}
\end{equation*}

It is evident

\noindent\textbf{Proposition 5}. {\it If the Krull dimension of an
algebra $A$ is $\varrho(A)$, then on any variety $V=\{a=0,$ $\
b<0\}$ with fixed $a,b\in A$} ($b$ {\it has no effect on the
mentioned variety}) {\it in the algebra $A$ there are not more than
$\varrho(A)$ algebraically independent elements $($possibly no more
than $\varrho(A)$ elements which form an integer algebraic basis$)$
of this algebra}.

\section{Hilbert series for Sibirsky graded algebras
$\boldsymbol{S_{1,m_1,m_2,...,m_{\ell}}}$ and
$\boldsymbol{SI_{1,m_1,m_2,...,m_{\ell}}}$}

\ \ \ \ According to Proposition 1 for the spaces of the algebra
$S_{1,m_1,m_2,...,m_{\ell}}$ from \eqref{eq501} we have
$dim_{\mathbb{R}}S_{1,m_1,m_2,...,m_{\ell}}^{(d)}<\infty$. Then,
following [7], by {\it the generalized Hilbert series} of the
algebra $S_{1,m_1,m_2,...,m_{\ell}}$ we mean a formal series
\begin{equation}
\begin{gathered}\label{eq12}
H(S_{1,m_1,m_2,...,m_{\ell}},u,z_0,z_1,...,z_{\ell})=\sum_{(d)}dim_{\mathbb{R}}S_{1,m_1,m_2,...,m_{\ell}}^{(d)}u^{\delta}z_0^{d_0}z_1^{d_1}...z_{\ell}^{d_{\ell}},
\end{gathered}
\end{equation}
which is said to reflect a $u,z$--graduation of the considered
algebra.

From the definition of the algebra of invariants
$SI_{1,m_1,m_2,...,m_{\ell}}$ and \eqref{eq12} it follows that
\begin{equation}
\begin{gathered}\label{eq13}
H(SI_{1,m_1,m_2,...,m_{\ell}},z_0,z_1,...,z_{\ell})=H(S_{1,m_1,m_2,...,m_{\ell}},0,z_0,z_1,...,z_{\ell}),
\end{gathered}
\end{equation}
and {\it the common Hilbert series} will be written respectively
\begin{equation}
\begin{gathered}\label{eq14}
H_{S_{1,m_1,m_2,...,m_{\ell}}}(u)=H(S_{1,m_1,m_2,...,m_{\ell}},u,u,u,...,u),\\
H_{SI_{1,m_1,m_2,...,m_{\ell}}}(z)=H(SI_{1,m_1,m_2,...,m_{\ell}},z,z,...,z).
\end{gathered}
\end{equation}
The last series contain meaningfull information about asymptotic
character of the behavior of the considered algebras.

The method of construction of the generalized Hilbert series
\eqref{eq12}--\eqref{eq14} for the algebras
$S_{1,m_1,m_2,...,m_{\ell}}$ and $SI_{1,m_1,m_2,...,m_{\ell}}$ was
developed in [7].

For example, the generalized Hilbert series for the algebras
$S_{0,1}$ and $SI_{0,1}$ of unimodular comitants and invariants  of
the system $s(0,1)$ have, respectively, the forms
\begin{equation*}
\begin{gathered}
H(S_{0,1},u,z_0,z_1)=\frac{1+uz_0z_1}{(1-uz_0)(1-z_1)(1-z_1^2)(1-z_0^2z_1)(1-u^2z_1)},\\
H(SI_{0,1},z_0,z_1)=\frac{1}{(1-z_1)(1-z_1^2)(1-z_0^2z_1)},
\end{gathered}
\end{equation*}
and the corresponding common Hilbert series will be written as
\begin{equation*}
\begin{gathered}
H_{S_{0,1}}(u)=\frac{1-u+u^2}{(1-u)^2(1-u^2)(1-u^3)^2},\
H_{SI_{0,1}}(z)=\frac{1}{(1-z)(1-z^2)(1-z^3)}.
\end{gathered}
\end{equation*}

\noindent\textit{Remark} 2.  We note, following $[21]$, that the
Krull dimension  $\varrho(S_{1,m_1,m_2,...,m_{\ell}})$ (respectively
$\varrho(SI_{1,m_1,m_2,...,m_{\ell}}))$ of the graded algebra
$S_{1,m_1,m_2,...,m_{\ell}}$ (respectively
$SI_{1,m_1,m_2,...,m_{\ell}})$ is equal to the multiplicity of the
pole of the common Hilbert series
$H_{S_{1,m_1,m_2,...,m_{\ell}}}(u)$ (respectively
$H_{SI_{1,m_1,m_2,...,m_{\ell}}}(z))$ at the unit.

For example, considering the above mentioned common Hilbert series
$H_{S_{0,1}}(u)$ and $H_{SI_{0,1}}(z)$ for the Krull dimension of
the algebras $S_{0,1}$ and $SI_{0,1}$ we obtain $\varrho(S_{0,1})=5$
and $\varrho(SI_{0,1})=3$, respectively.

In other cases, when there is no explicit form of the common Hilbert
series, but the power series expansion is known, then we can use the
following

\noindent\textit{Remark} 3.  Accept that the comparison of series
with non-negative coefficients is performed coefficient-wise $(\sum
a_nt^n\leq\sum b_nt^n\Leftrightarrow a_n\leq b_n;\ \forall n)$.
Taking this into account, if for commutative graded algebras $A$ and
$B$ we have
\begin{equation}
\begin{gathered}\label{eq15}
H_{A}(t)\leq H_{B}(t),
\end{gathered}
\end{equation}
then for their Krull dimensions  we also have
$\varrho(A)\leq\varrho(B)$.

It is also evident that if for the common Hilbert series of a
commutative graded algebra $A$ we have
\begin{equation}
\begin{gathered}\label{eq16}
H_{A}(t)\leq \frac{C}{(1-t)^m},
\end{gathered}
\end{equation}
where $C$ is a fixed constant, then we obtain $\varrho(A)\leq m$.

The extended theory and bibliography about Hilbert series for graded
algebras can be found in [22].

\section{Lie algebras of operators admitted by polynomial dif-ferential systems}

\ \ \ \ It is shown in [7] that any differential system
$s(m_0,m_1,m_2,...,m_{\ell})$ from \eqref{eq1} admits a
four-dimensional reductive Lie algebra $L_4$, which consists of
operators
\begin{equation}
\begin{gathered}\label{eq17}
X_1=x\frac{\partial}{\partial x}+D_1,\ X_2=y\frac{\partial}{\partial
x}+D_2,\ X_3=x\frac{\partial}{\partial y}+D_3,\
X_4=y\frac{\partial}{\partial y}+D_4,
\end{gathered}
\end{equation}
where the differential operators $D_1,D_2,D_3,D_4$ are operators of
the representation of the center-affine group $GL(2,\mathbb{R})$ in
the space of the coefficients of the polynomials $P_{m_i}$ and
$Q_{m_i}\ (i=\overline{1,\ell})$ of the system \eqref{eq1}.

In [7] it is proved

\noindent\textbf{Theorem 1}. {\it For a polynomial $k$ in the
coefficients of the system $s(m_0,m_1,m_2,...,m_{\ell})$ from
\eqref{eq1} and phase variables $x,y$ to be a center-affine comitant
of this system with the weight $g$, it is necessary and sufficient
that it satisfies the equations}
\begin{equation*}
\begin{gathered}
X_1(k)=X_4(k)=-gk,\ X_2(k)=X_3(k)=0.
\end{gathered}
\end{equation*}

With the help of this theorem and properties of rational absolute
center-affine comitants of the system \eqref{eq1} from [7],
following the classical theory of these invariants [16], it can be
shown that for the number of elements in an algebraic basis of
center-affine comitants of the system $s(m_0,m_1,m_2,...,m_{\ell})$
the following formula holds:
\begin{equation}
\begin{gathered}\label{eq18}
\varrho=2\left(\sum_{i=0}^{\ell}m_i+\ell\right)+1.
\end{gathered}
\end{equation}

In the theory of center-affine comitants of polynomial differential
systems [6] it is shown that if $S$ is a {\it semi-invariant} in the
center-affine comitant $k$, then
\begin{equation}
\begin{gathered}\label{eq19}
k=Sx^{\delta}-D_3(S)x^{\delta-1}y+\frac{1}{2!}D_3^2(S)x^{\delta-2}y^2+...+\frac{(-1)^{\delta}}{\delta!}D_3^{\delta}(S)y^{\delta},
\end{gathered}
\end{equation}
where $D_3$ is defined in [7].

\noindent\textit{Remark} 4. [6] With the help of this equality it
can be shown that the center-affine comitants
$k_1,k_2,...,k_{\varrho(S_{m_0,m_1,m_2,...,m_{\ell}})}\in
S_{m_0,m_1,m_2,...,m_{\ell}}$ which belong to the system
$s(m_0,m_1,$ $m_2,...,m_{\ell})$ are algebraically independent if
and only if their semi-invariants are algebraically independent.

\section{An invariant variety in the center-focus problem of the system
$\boldsymbol{s(1,m_1,m_2,...,m_{\ell})}$}

\ \ \ \ The center-focus problem for systems of the form \eqref{eq1}
has the following classical formulation: for an infinite system of
polynomials
\begin{equation}
\begin{gathered}\label{eq20}
\{(x^2+y^2)^k\}_{k=1}^{\infty}
\end{gathered}
\end{equation}
there exists a function
\begin{equation}
\begin{gathered}\label{eq21}
U(x,y)=x^2+y^2+\sum_{k=3}^{\infty}f_k(x,y),
\end{gathered}
\end{equation}
where $f_{k}(x,y)$ are homogeneous polynomials of degree $k$ in
$x,y$, and such constants
\begin{equation*}
\begin{gathered}
L_1,L_2,...,L_k,...
\end{gathered}\eqno(2)
\end{equation*}
that the identity
\begin{equation}
\begin{gathered}\label{eq22}
\frac{dU}{dt}=\sum_{k=1}^{\infty}L_{k}(x^2+y^2)^{k+1}
\end{gathered}
\end{equation}
 (with respect to $x$ and $y$) holds along the trajectories of the system
\begin{equation}
\begin{gathered}\label{eq23}
\dot{x}=y+\sum_{i=1}^{\ell}P_{m_i}(x,y),\
\dot{y}=-x+\sum_{i=1}^{\ell}Q_{m_i}(x,y).
\end{gathered}
\end{equation}
The constants \eqref{eq2} are polynomials in coefficients of the
system \eqref{eq23}, and are called {\it focal quantities}.

We note that the algebra $S_{1,m_1,m_2,...,m_{\ell}}$ for any
differential system $s(1,m_1,m_2,...,m_{\ell})$, written in the form
\begin{equation}
\begin{gathered}\label{eq24}
\dot{x}=cx+dy+\sum_{i=1}^{\ell}P_{m_i}(x,y),\
\dot{y}=ex+fy+\sum_{i=1}^{\ell}Q_{m_i}(x,y)
\end{gathered}
\end{equation}
contains among its generators the polynomials
\begin{equation}
\begin{gathered}\label{eq25}
i_1=c+f,\ i_2=c^2+2de+f^2,\ k_2=-ex^2+(c-f)xy+dy^2,
\end{gathered}
\end{equation}
which are given already in \eqref{eq10}.

\noindent\textit{Remark} 5. We note that the set
\begin{equation}
\begin{gathered}\label{eq26}
\mathcal{V}=\{i_1=c+f=0,\ Discr(k_2)=2i_2-i_1^2<0\}
\end{gathered}
\end{equation}
is a Sibirsky invariant variety for center and focus for the system
\eqref{eq24}, because the comitant $k_2$ from \eqref{eq25} through a
real center-affine transformation of the plane $xOy$ can be brought
to the form
\begin{equation}
\begin{gathered}\label{eq27}
x^2+y^2,
\end{gathered}
\end{equation}
and the system \eqref{eq24} can be brought to the form \eqref{eq23}
$[5]$, for which the roots of the characteristic equation are
imaginary, i.\,e. the origin of coordinates for this system is a
singular point of the second type (center or focus).

Considering Remark 5 we have

\noindent\textit{Remark} 6. Taking into account the comitant $k_2$
from \eqref{eq25} and the fact that its expression through a real
center-affine transformation on the invariant variety $\mathcal{V}$
can be brought to the form \eqref{eq27}, then formally this variety
for the system \eqref{eq24} can be written as
\begin{equation}
\begin{gathered}\label{eq28}
\mathcal{V}=\{f=-c\}\cup\{c=0,\ d=-e=1\}.
\end{gathered}
\end{equation}

\section{Null focal pseudo-quantity of the system (25) and
relations between the quantities $\boldsymbol{G_k}$ and the focal
quantities $\boldsymbol{L_k}$ of the system (24)}

\ \ \ \ Let us consider for the system \eqref{eq24} the identity
\begin{equation}
\begin{gathered}\label{eq29}
\left[cx+dy+\sum_{i=1}^{\ell}P_{m_i}(x,y)\right]\frac{\partial{U}}{\partial{x}}+\left[ex+fy+\sum_{i=1}^{\ell}Q_{m_i}(x,y)\right]\frac{\partial{U}}
{\partial{y}}=\sum_{k=1}^{\infty}G_{k}k_{2}^{k+1},
\end{gathered}
\end{equation}
where
\begin{equation}
\begin{gathered}\label{eq30}
U(x,y)=k_2+\sum_{r=3}^{\infty}F_r(x,y),
\end{gathered}
\end{equation}
$(k_2\not\equiv0$ from \eqref{eq25}), which splits by powers of $x$
and $y$ into an infinite number of algebraic equations, where the
variables are the coefficients of the homogeneous polynomials
$F_r(x,y)$ of degree $r$ in $x,y$, and also the quantities
$G_1,G_2,...,G_k,...$.

For any system \eqref{eq24} from the identity \eqref{eq29} with
$k_2$ from \eqref{eq25} we find that the first three equations have
the following form:
\begin{equation*}
\begin{gathered}
x^2:\ e(c+f)=0,\ xy:\ (c-f)(c+f)=0,\ y^2:\ d(c+f)=0.
\end{gathered}
\end{equation*}
These equalities are equivalent to one of two sets of the
conditions: $1)\ c+f=0$;  $2)\ e=c-f=$ $=d=0$. Since
$k_2\not\equiv0$, then, according to \eqref{eq25}, these conditions
are equivalent to the condition $c+f=0$, which is contained in the
variety $\mathcal{V}$ from \eqref{eq26}.

In this way from Remark 5 (6) and formulation of the center-focus
problem for the system \eqref{eq23} we conclude: for $L_k$ from
\eqref{eq2} and $G_k$ from \eqref{eq29} the following equalities
take place:
\begin{equation}
\begin{gathered}\label{eq31}
L_k=G_k|_{\mathcal{V}}\ (k=1,2,3,...),
\end{gathered}
\end{equation}
where $\mathcal{V}$ is from \eqref{eq26}.

Hereafter some concretizations for these equalities will be done.

From the above mentioned follows

\noindent\textit{Remark} 7. The identity \eqref{eq29} with function
\eqref{eq30} on the variety $\mathcal{V}$ from \eqref{eq26}
guarantees that the system \eqref{eq24} has at the origin of
coordinates a singular point of the second type (center or focus).

We denote the expression $c+f$, which is contained in the variety
$\mathcal{V}$ from \eqref{eq26}, by
\begin{equation}
\begin{gathered}\label{eq32}
G_0\equiv i_1=c+f,
\end{gathered}
\end{equation}
and will call it {\it the null focal pseudo-quantity}. We note that
$G_0$ from \eqref{eq32} is a center-affine (unimodular) invariant of
the system $s(1,m_1,m_2,...,m_{\ell})$ of the type
\begin{equation*}
\begin{gathered}
(0,1,{\underbrace{0,...,0}).}\atop{\ \ \ \ \  \ell}
\end{gathered}
\end{equation*}
To get a more clear idea about the quantities $G_1,G_2,...,G_k,...$
from the identity \eqref{eq29} with the function \eqref{eq30}, we
write the remaining equations, in which this identity is splitted by
powers $x^3,x^2y,xy^2,y^3,...$ without taking into consideration the
equality $i_1=c+f=0$ on the variety $\mathcal{V}$.

To explain the further way of implementation of this scenario, we
consider the identity \eqref{eq29} with unknown constants
$G_1,G_2,...$ for the example of the simplest differential system
$s(1,2)$ with the quadratic nonlinearities
\begin{equation}
\begin{gathered}\label{eq33}
\dot{x}=cx+dy+gx^2+2hxy+ky^2,\\
\dot{y}=ex+fy+lx^2+2mxy+ny^2,
\end{gathered}
\end{equation}
with the finitely defined graded algebra of unimodular comitants
$S_{1,2}\ [7]$. For this algebra we write the function \eqref{eq30}
as
\begin{equation}
\begin{gathered}\label{eq34}
U(x,y)=k_2+a_0x^3+3a_1x^2y+3a_2xy^2+a_3y^3+b_0x^4+4b_1x^3y+\\+6b_2x^2y^2+4b_3xy^3+b_4y^4+c_0x^5+5c_1x^4y+10c_2x^3y^2+\\+10c_3x^2y^3+
5c_4xy^4+c_5y^5+d_0x^6+6d_1x^5y+15d_2x^4y^2+\\+20d_3x^3y^3+
15d_4x^2y^4+6d_5xy^5+d_6y^6+e_0x^7+7e_1x^6y+\\+21e_2x^5y^2+35e_3x^4y^3+21e_5x^2y^5+7e_6xy^6+e_7y^7+f_0x^8+\\+8f_1x^7y+
28f_2x^6y^2+56f_3x^5y^3+70f_4y^4+56f_5x^3y^5+\\+28f_6x^2y^6+8f_7xy^7+f_8y^8+...,
\end{gathered}
\end{equation}
where $k_2$ is from \eqref{eq25} and $a_0,a_1,...,f_7,f_8,...$ are
unknown constants.  Then without taking into consideration the
variety $\mathcal{V}$, the identity \eqref{eq29} along the
trajectories of the system \eqref{eq33} with the function
\eqref{eq34} splits into the following systems of equations
\begin{equation}
\begin{split}\label{eq35}
x^3:& \  3ca_0+3ea_1=2eg-(c-f)l,
\\ x^2y:& \
3da_0+3(2c+f)a_1+6ea_2=(f-c)(g+2m)-2dl+4eh,
\\ xy^2:& \  6da_1+3(2f+c)a_2+3ea_3=(f-c)(2h+n)+2ek-4dm,
\\ y^3: & \ 3da_2+3fa_3=(f-c)k-2dn;
\end{split}
\end{equation}
\begin{equation}
\begin{split}\label{eq36}
x^4:&  \  4cb_0+4eb_1-e^2G_1=-3ga_0-3la_1,
\\  x^3y:&  \ 4db_0+4(f+3c)b_1+12eb_2+2e(c-f)G_1=-6ha_0-\\ &-6(g+m)a_1-6la_2,
\\  x^2y^2:& \  12db_1+12(c+f)b_2+12eb_3+[2de-(c-f)^2]G_1=\\ &=-3ka_0-3(4h+n)a_1-
3(g+4m)a_2-3la_3,
\\ xy^3:& \
12bd_2+4(3f+c)b_3+4eb_4+2d(f-g)G_1=-6ka_1-\\ &-6(h+n)a_2-6ma_3,
\\ y^4:&  \ 4db_3+4fb_4-d^2G_1=-3ka_2-3na_3;
\end{split}
\end{equation}
\begin{equation*}
\begin{split}
x^5:&  \  5cc_0+5ec_1=-4gb_0-4lb_1,
\\ x^4y:&  \ 5dc_0+5(4c+f)c_1+20ec_2=-8hb_0-4(3g+2m)b_1-\\ &-12lb_2,
\end{split}
\end{equation*}
\begin{equation}
\begin{split}\label{eq37}
x^3y^2:& \
20dc_1+10(3c+2f)c_2+30ec_3=-4kb_0-4(6h+n)b_1-\\
&-12(g+2m)b_2-12lb_3,
\\ x^2y^3:& \ 30dc_2+10(2c+3f)c_3+20ec_4=-12kb_1-12(2h+\\ &+n)b_2-4(g+6m)b_3-4lb_4,
\\xy^4:& \ 20dc_3+5(c+4f)c_4+5ec_5=-12kb_2-4(2h+3n)b_3-\\ &-8mb_4,
\\ y^5:&  \ 5dc_4+5fc_5=-4kb_3-4nb_4;
\end{split}
\end{equation}
\begin{equation}
\begin{split}\label{eq38}
x^6:&  \  6cd_0+6ed_1+e^3G_2=-5gc_0-5lc_1,
\\ x^5y:&  \ 6dd_0+6(5c+f)d_1+30ed_2+3e^2(f-c)G_2=-10hc_0-\\ &-10(2g+m)c_1-20lc_2,
\\ x^4y^2:& \
30dd_1+30(2c+f)d_2+60ed_3+3e[(c-f)^2-de]G_2=\\
&=-5kc_0-5(8h+n)c_1-10(3g+4m)c_2-30lc_3,
\\x^3y^3:& \ 60dd_2+60(c+f)d_3+60ed_4+(f-c)[(c-f)^2-\\
&-6de]G_2=-20kc_1-20(3h+n)c_2-20(g+3m)c_3-20lc_4,
\\  x^2y^4:& \ 60dd_3+30(c+2f)d_4+30ed_5+3d[de-(c-f)^2]G_2=\\ &=-30kc_2-10(4h+3n)c_3-5(g+8m)c_4-5lc_5,
\\ xy^5:& \
30dd_4+6(c+5f)d_5+6ed_6+3d^2(f-c)G_2=-20kc_3-\\
&-10(h+2n)c_4-10mc_5,
\\ y^6:&  \ 6dd_5+6fd_6-d^3G_2=-5kc_4-5nc_5;
\end{split}
\end{equation}
\begin{equation}
\begin{split}\label{eq39}
x^7:&  \  7ce_0+7ee_1=-6gd_0-6ld_1,
\\ x^6y:&  \ 7de_0+7(6c+f)e_1+42ee_2=-12hd_0-6(5g+2m)d_1-\\ &-30ld_2,
\\ x^5y^2:& \ 42de_1+7(15c+6f)e_2+105ee_3=-6kd_0-6(10h+\\ &+n)d_1-60(g+m)d_2-60ld_3,
\\ x^4y^3:& \ 105de_2+5(28c+21f)e_3+140ee_4=-30kd_1-30(4h+\\ &+n)d_2-60(g+2m)d_3-60ld_4,
\\ x^3y^4:& \ 140de_3+35(3c+4f)e_4+105ee_5=-60kd_2-60(2h+\\ &+n)d_3-30(g+4m)d_4-30ld_5,
\\ x^2y^5:& \ 105de_4+7(6c+15f)e_5+42ee_6=-60kd_3-60(h+\\ &+n)d_4-6(g+10m)d_5-6ld_6,
\\xy^6:& \ 42de_5+7(c+6f)e_6+7ee_7=-30kd_4-6(2h+5n)d_5-\\ &-12md_6,
\\ y^7:&  \ 7de_6+7fe_7=-6kd_5-6nd_6;
\end{split}
\end{equation}
\begin{equation*}
\begin{split}
x^8:&  \  8cf_0+8ef_1-e^4G_3=-7ge_0-7le_1,
\\ x^7y:&  \ 8df_0+8(7c+f)f_1+56ef_2+4e^3(c-f)G_3=\\
&=-14he_0-14(3g+m)e_1-42le_2,
\end{split}
\end{equation*}
\begin{equation}
\begin{split}\label{eq40}
x^6y^2:& \ 56df_1+56(3c+f)f_2+168ef_3+2e^2[2de-3(c-f)^2]G_3=\\
&=-7ke_0-7(12h+n)e_1-21(5g+4m)e_2-105le_3,
\\ x^5y^3:& \ 168df_2+56(5c+3f)f_3+280ef_4+4e(f-c)[3de-(c-\\ &-f)^2]G_3=-42ke_1-42(5h+n)e_2-70(2g+3m)e_3-140le_4,
\\ x^4y^4:& \ 280df_3+280(c+f)f_4+280ef_5+[12de(c-f)^2-6d^2e^2-\\ &-(c-f)^4]G_3=-105ke_2-35(8h+3n)e_3-35(3g+\\ &+8m)e_4-105le_5,
\\ x^3y^5:& \ 280df_4+56(3c+5f)f_5+168ef_6+4d(f-c)[(c-f)^2-\\ &-3de]G_3=-140ke_3-70(3h+2n)e_4-42(g+5m)e_5-42le_6,
\\ x^2y^6:& \ 168df_5+56(c+3f)f_6+56ef_7+2d^2[2de-3(c-f)^2]G_3=\\ &=-105ke_4-21(4h+5n)e_5-7(g+12m)e_6-7le_7,
\\ xy^7:& \ 56df_6+8(c+7f)f_7+8ef_8+4d^3(f-c)G_3=-42ke_5-
\\  &-14(h+3n)e_6-14me_7,
\\ y^8:&  \ 78df_7+8ff_8-d^4G_3=-7ke_6-7ne_7.
\end{split}
\end{equation}

It is evident that the linear systems of equations
\eqref{eq35}--\eqref{eq40} in variables
$a_0,a_1,a_2,a_3,b_0,b_1,...,b_4,$
$c_0,c_1,...,c_5,d_0,d_1,...,d_6,e_0,e_1,...,e_7,f_0,f_1,...,f_8,...,G_1,G_2,G_3,...$
can be extended by adding, after the last equation from
\eqref{eq40}, an infinite number of equations, obtained from the
equality of coefficients of $x^{\alpha} y^{\beta}$ for
$\alpha+\beta>8$ in the identity \eqref{eq29} along the trajectories
of the system \eqref{eq33}.

\section{Determining the quantities $\boldsymbol{G_1,G_2,G_3}$ from the systems
(36)--(41) and the corresponding focal quantities}

\ \ \ \ To obtain the quantity  $G_1$ we write the equations
\eqref{eq35}--\eqref{eq36} in the matrix form
\begin{equation}
\begin{gathered}\label{eq41}
A_{1} B_{1}=C_{1},
\end{gathered}
\end{equation}
where
{\footnotesize
$$
\begin{gathered}
A_{1}\!=\!\left(\!\!
\begin{array}{cccccccccc}
  3c&3e&0&0&0&0&0&0&0&0\\
  3d&3(2c\!+\!f)&6e&0&0&0&0&0&0&0\\
  0&6d&3(2c\!+\!f)&3e&0&0&0&0&0&0\\
  0&0&3d&3f&0&0&0&0&0&0\\
  3g&3l&0&0&4c&4e&0&0&0&-e^2\\
  6h&6(g\!+\!m)&6l&0&4d&4(f\!+\!3c)&12e&0&0&2e(c\!-\!f)\\
  3k&3(4h\!+\!n)&3(g\!+\!4m)&3l&0&12d&12(c\!+\!f)&12e&0&2de\!-\!(c\!\!f)^2\\
  0&6k&6(h\!+\!n)&6m&0&0&12d&4(3f\!+\!c)&4l&2d(f\!-\!c)\\
  0&0&3k&3n&0&0&0&4d&4f&-d^2
\end{array}\!\!\!
\right),\\
\end{gathered}
$$
\begin{equation}\label{eq42}
\begin{gathered}
B_{1}=\left(%
\begin{array}{c}
  a_0\\
  a_1\\
  a_2\\
  a_3\\
  b_0\\
  b_1\\
  b_2\\
  b_3\\
  b_4\\
  G_1
\end{array}%
\right),\quad
C_{1}=\left(%
\begin{array}{c}
  2eg+(f-c)l \\
  (f-c)(g+2m)-2dl+4eh \\
  (f-c)(2h+n)+3ek-4dm\\
  (f-c)k-2dn\\
  0\\
  0\\
  0\\
  0\\
  0\\
\end{array}%
\right).
\end{gathered}
\end{equation}
}

Sice the dimension of the matrix $A_1$ is $9\times10$, clearly we
have at least one free parameter. Therefore choosing as a free
parameter $b_i\ (i\in\{0,1,...,4\})$ with the help of Cramer's rule
from the system \eqref{eq41} for each fixed $i$ we obtain
\begin{equation}
\begin{gathered}\label{eq43}
G_1=\frac{G_{1,i}+B_{1,i}b_i}{\sigma_{1,i}}
\end{gathered}
\end{equation}
where $G_{1,i},B_{1,i},\sigma_{1,i}$ are polynomials in the
coefficients of the system \eqref{eq33}, and $b_i$ are undetermined
coefficients of the function $U(x,y)$ from \eqref{eq34}.

By studying the matrices \eqref{eq42} of the system \eqref{eq41} we
conclude that $G_{1,i}$ from \eqref{eq43} are homogeneous
polynomials of degree $8$ with respect to the linear part, and of
degree $2$ with respect to the quadratic part of the system
\eqref{eq33}.

Because $G_{1,i}$ from \eqref{eq43} are homogeneous polynomials in
the coefficients of the system \eqref{eq33}, then, according to
[6,23], for $i=0,1,2,3,4$ we can determine respectively and
isobarity
\begin{equation*}
\begin{gathered}
(3,-1),\ (2,0),\ (1,1),\ (0,2),\ (-1,3).
\end{gathered}
\end{equation*}
According to the formula \eqref{eq501} (for the system \eqref{eq33}
and the theory of invariants of differential systems [5,6]) it
suggests that the numerators of the fractions \eqref{eq43} can be
coefficients in comitants of the weight $-1$ of the type $(4,8,2)$. This
means that according to \eqref{eq19} with the help of the Lie
differential operator $D_3$ for the system \eqref{eq33} from [7] and
the numerator of the fraction \eqref{eq43} we obtain a redefined
system of four linear non-homogeneous differential equations
\begin{equation}
\begin{gathered}\label{eq44}
D_3(G_{1,0}+B_{1,0}b_0)=G_{1,1}+B_{1,1}b_1,\
D_3(G_{1,1}+B_{1,1}b_1)=-G_{1,2}-B_{1,2}b_2,\\
-D_3(G_{1,2}+B_{1,2}b_2)=G_{1,3}+B_{1,3}b_3,\
D_3(G_{1,3}+B_{1,3}b_3)=-G_{1,4}-B_{1,4}b_4
\end{gathered}
\end{equation}
with five unknowns $b_0,b_1,b_2,b_3,b_4$. We can note that a
particular solution to this system is $b_0=b_1=b_2=b_3=b_4=0$, for
which the polynomial
\begin{equation}
\begin{gathered}\label{eq45}
f_4^{\prime}(x,y)=G_{1,0}x^4+4G_{1,1}x^3y+2G_{1,2}x^2y^2+4G_{1,3}xy^3+G_{1,4}y^4
\end{gathered}
\end{equation}
is a center-affine comitant of the system \eqref{eq33}. This fact is
also confirmed by Theorem 1 with the operators $X_1-X_4$ from [7]
for the system \eqref{eq33}, for which
\begin{equation*}
\begin{gathered}
X_1(f_4^{\prime})=X_4(f_4^{\prime})=f_4^{\prime},\
X_2(f_4^{\prime})=X_3(f_4^{\prime})=0.
\end{gathered}
\end{equation*}

Similarly, one can see that another particular solution for the
system \eqref{eq44} is given by the following expressions:
\begin{equation*}
\begin{split}
& b_0 =\frac{-e (g^2 + 2 h l + m^2)}{3 c^2 - 4 d e + 10 c f + 3 f^2},\\
& b_1 = \frac{(c-f) (g^2 + 2 h l + m^2)-
  2 e (g h + k l + h m  + m n)}{4 (3 c^2 - 4 d e + 10 c f + 3 f^2)},\\
& b_2 = \frac{2(c-f)(g h + k l + h m + m n) - e(h^2 + 2 k m + n^2)+
d (g^2 + 2 h l +
m^2)}{6 (3 c^2 - 4 d e + 10 c f + 3 f^2)},\\
& b_3 = \frac{(c - f) (h^2 + 2 k m + n^2) + 2 d (g h + k l + h m  +
m n)}{
 4 (3 c^2 - 4 d e + 10 c f + 3 f^2)},\\
& b_4 = \frac{d (h^2 + 2 k m + n^2)}{3 c^2 - 4 d e + 10 c f + 3
f^2},
\end{split}
\end{equation*}
whose denominators are different from zero on the variety
$\mathcal{V}$ from \eqref{eq26}. They define the center-affine
comitant
\begin{equation}
\begin{gathered}\label{eq46}
f_4^{\prime\prime}(x,y)=(G_{1,0}+B_{1,0}b_0)x^4+4(G_{1,1}+B_{1,1}b_1)x^3y+2(G_{1,2}+\\+B_{1,2}b_2)x^2y^2+4(G_{1,3}+B_{1,3}b_3)xy^3+(G_{1,4}+B_{1,4}b_4)y^4.
\end{gathered}
\end{equation}

It is evident that the differential system \eqref{eq44} has an
infinite number of solutions $b_0,b_1,b_2,b_3,b_4$, which define
center-affine comitants of the type \eqref{eq46}.

In view of the above, the comitants \eqref{eq45}--\eqref{eq46}
belong to the space
\begin{equation*}
\begin{gathered}
S_{1,2}^{(4,8,2)}.
\end{gathered}
\end{equation*}
Remark that the comitants \eqref{eq45}--\eqref{eq46} on the variety
$\mathcal{V}$ from \eqref{eq26} for the system \eqref{eq33} have the
following form:
\begin{equation}
\begin{gathered}\label{eq47}
f_4^{\prime}(x,y)|_{\mathcal{V}}=f_4^{\prime\prime}(x,y)|_{\mathcal{V}}=-8L_1(x^2+y^2)^2,\
\end{gathered}
\end{equation}
where
\begin{equation*}
\begin{gathered}
L_1=\frac{1}{2}\left[g(l-h)-k(h+n)+m(l+n)\right]
\end{gathered}
\end{equation*}
is the first focal quantity of the system \eqref{eq33} on the
invariant variety $\mathcal{V}$ (see [4, p.~110]).

Similarly to the previous case, for determining the quantity $G_{2}$
we write the equations \eqref{eq35}--\eqref{eq38} in the matrix form
\begin{equation}
\begin{gathered}\label{eq48}
A_{2} B_{2}=C_{2},
\end{gathered}
\end{equation}
from which we find
\begin{equation}
\begin{gathered}\label{eq49}
G_2=\frac{G_{2,i,j}+B_{2,i,j}b_i+D_{2,i,j}d_j}{\sigma_{2,i,j}},\
(i=\overline{0,4},\ j=\overline{0,6}).
\end{gathered}
\end{equation}
By studying the matrix equality \eqref{eq48} we obtain that $deg
G_{2,i,j}=24$, and using the system \eqref{eq35}--\eqref{eq38} we
obtain that $G_{2,i,j}$ from \eqref{eq49} has the type $(0,20,4)$,
i.\,e. $G_{2,i,j}$ are homogeneous polynomials of degree $20$ in
coefficients of the linear part and of degree $4$ in coefficients of
the quadratic part of the system $s(1,2)$ from \eqref{eq33}.

Computing the expressions $G_{2,i,j}$ for each $i=\overline{0,4}$
and $j=\overline{0,6}$, according to [6,23], we obtain for their
isobarity the following table:

\ \ \ \ \ \ \ \ \ \ \ \ \ \ \ \ \ \ \ \ \ \ \ \ \ \ \ \ \ \ \ \ \ \
\ \ \ \ \ \ \ \ \ \ \ \ \ \ \ \ \ \ \ \ \ \ \ \ \ \ \ \ \ \ \ \ \ \
\ \ \ \ Table 1
\begin{equation*}
\begin{gathered}
\begin{tabular}{|c|c|c|c|c|c|c|c|}
  \hline
  $G_{2,i,j}$ & $d_0$ & $d_1$ & $d_2$ & $d_3$ & $d_4$ & $d_5$ & $d_6$\\
  \hline
  $b_0$ & (7,-3) & (6,-2) & (5,-1) & (4,0) & (3,1) & (2,2) & (1,3)\\ \hline
  $b_1$ & (6,-2) & (5,-1) & (4,0) & (3,1) & (2,2) & (1,3) & (0,4)\\ \hline
  $b_2$ & (5,-1) & (4,0) & (3,1) & (2,2) & (1,3) & (0,4) & (-1,5)\\ \hline
  $b_3$ & (4,0) & (3,1) & (2,2) & (1,3) & (0,4) & (-1,5) & (-2,6)\\ \hline
  $b_4$ & (3,1) & (2,2) & (1,3) & (0,4) & (-1,5) & (-2,6) & (-3,7)\\ \hline
\end{tabular}
\end{gathered}
\end{equation*}
By studying the isobarity of $G_{2,i,j}$ top-down for each line of
this table, according to the theory of invariants of differential
systems [5,6], we find that the numerators of the fraction
\eqref{eq49} can be coefficients in center-affine comitants with the
corresponding weights $-3,-2,-1,0,1$. Using these weights and the
formula \eqref{eq501} for the system \eqref{eq33}, as well as the
fact that $G_{2,i,j}$ have the type $(0,20,4)$, we obtain that the
mentioned comitants correspond to the types
\begin{equation}
\begin{gathered}\label{eq50}
(10,20,4),\ (8,20,4),\ (6,20,4),\ (4,20,4),\ (2,20,4).
\end{gathered}
\end{equation}

As the quantity $G_{2}$ in \eqref{eq29} is the coefficient in front
of the homogeneity of degree $6$ in the phase variables $x$ and $y$,
then it is logical to choose from \eqref{eq50} the type
\begin{equation}
\begin{gathered}\label{eq51}
(6,20,4),
\end{gathered}
\end{equation}
which corresponds to the expression $G_{2,2,j} \
 (j=\overline{0,6})$ in Table 1.

This means that according to \eqref{eq19} using the Lie differential
operator $D_3$ for the system \eqref{eq33} from [7] and the
numerator of the fraction \eqref{eq49} for fixed $i=2$, we obtain
one redefined system of six linear non-homogeneous differential
equations
\begin{equation}
\begin{gathered}\label{eq52}
D_3(G_{2,2,0}+B_{2,2,0}b_0+D_{2,2,0}d_0)=-(G_{2,2,1}+B_{2,2,1}b_1+D_{2,2,1}d_1),\\
-D_3(G_{2,2,1}+B_{2,2,1}b_1+D_{2,2,1}d_1)=G_{2,2,2}+B_{2,2,2}b_2+D_{2,2,2}d_2,\\
D_3(G_{2,2,2}+B_{2,2,2}b_2+D_{2,2,2}d_2)=-(G_{2,2,3}+B_{2,2,3}b_3+D_{2,2,3}d_3),\\
-D_3(G_{2,2,3}+B_{2,2,3}b_3+D_{2,2,3}d_3)=G_{2,2,4}+B_{2,2,4}b_4+D_{2,2,4}d_4,\\
D_3(G_{2,2,4}+B_{2,2,4}b_4+D_{2,2,4}d_4)=-(G_{2,2,5}+B_{2,2,5}b_5+D_{2,2,5}d_5),\\
-D_3(G_{2,2,5}+B_{2,2,5}b_5+D_{2,2,5}d_5)=G_{2,2,6}+B_{2,2,6}b_6+D_{2,2,6}d_6,
\end{gathered}
\end{equation}
with eight unknowns $b_2,d_0,d_1,...,d_6$. From these six equations
it results that the expressions contained in them can be
coefficients in comitants of the type $(6,20,4)$. Observe that
obtaining an explicit form for solutions of the system \eqref{eq52}
is a difficult task. We will show the importance of homogeneities of
$G_{2,2,j}$ from \eqref{eq49} in obtaining the focal quantities for
the system \eqref{eq33} on the invariant variety $\mathcal{V}$  for
center and focus  from \eqref{eq26}. According to \eqref{eq51} the
system \eqref{eq52} defines center-affine comitants belonging to the
space
\begin{equation}
\begin{gathered}\label{eq53}
S_{1,2}^{(6,20,4)}.
\end{gathered}
\end{equation}
According to \eqref{eq19} and \eqref{eq52} such a comitant,
belonging to this space, can be written as
\begin{equation*}
\begin{gathered}
f_6^{\prime}(x,y)=(G_{2,2,0}+B_{2,2,0}b_2+D_{2,2,0}d_0)x^6-(G_{2,2,1}+B_{2,2,1}b_2+D_{2,2,1}d_1)x^5y+\\+\frac{1}{2!}(G_{2,2,2}+B_{2,2,2}b_2+D_{2,2,2}d_2)x^4y^2
-\frac{1}{3!}(G_{2,2,3}+B_{2,2,3}b_2+D_{2,2,3}d_3)x^3y^3+\\+\frac{1}{4!}(G_{2,2,4}+B_{2,2,4}b_2+D_{2,2,4}d_4)x^2y^4
-\frac{1}{5!}(G_{2,2,5}+B_{2,2,5}b_2+D_{2,2,5}d_5)xy^5+\\+\frac{1}{6!}(G_{2,2,6}+B_{2,2,6}b_2+D_{2,2,6}d_6)y^6.
\end{gathered}
\end{equation*}
Observe that on the variety $\mathcal{V}$ from \eqref{eq26} for the
system \eqref{eq33} the expressions $G_{2,2,j}\ (j=\overline{0,6})$
have the following expressions:
\begin{equation}
\begin{gathered}\label{eq54}
G_{2,2,0}|_{\mathcal{V}}=G_{2,2,2}|_{\mathcal{V}}=G_{2,2,4}|_{\mathcal{V}}=G_{2,2,6}|_{\mathcal{V}}=-2304L_2,\\
G_{2,2,1}|_{\mathcal{V}}=G_{2,2,3}|_{\mathcal{V}}=G_{2,2,5}|_{\mathcal{V}}=0
\end{gathered}
\end{equation}
where
\begin{equation*}
\begin{gathered}
24L_2=62 g^3 h - 2 g h^3 + 95 g^2 h k - 2 h^3 k + 38 g h k^2 + 5 h
k^3 -
 62 g^3 l +\\+ 27 g h^2 l - 39 g^2 k l + 29 h^2 k l - 15 g k^2 l -
 8 g h l^2 + 15 h k l^2 - 5 g l^3 +\\+ 53 g^2 h m + 66 g h k m +
 13 h k^2 m - 127 g^2 l m - 6 h^2 l m - 68 g k l m -\\- 15 k^2 l m -
 13 h l^2 m - 5 l^3 m + 6 g h m^2 + 6 h k m^2 - 63 g l m^2 -
 29 k l m^2 +\\+ 2 l m^3 + 6 g^3 n + 61 g h^2 n + 72 g^2 k n +
 63 h^2 k n + 33 g k^2 n + 5 k^3 n -\\- 10 g h l n + 68 h k l n -
 33 g l^2 n + 15 k l^2 n - 72 g^2 m n - 6 h^2 m n +\\+ 10 g k m n +
 8 k^2 m n - 66 h l m n - 38 l^2 m n - 61 g m^2 n - 27 k m^2 n +\\+
 2 m^3 n +72 g h n^2 + 127 h k n^2 - 72 g l n^2 + 39 k l n^2 -
 53 h m n^2 - \\-95 l m n^2 - 6 g n^3 + 62 k n^3 - 62 m n^3
\end{gathered}
\end{equation*}
is the second focal quantity of the system \eqref{eq33} on the
invariant variety $\mathcal{V}$ for center and focus (see [4, p.~110]).

Now we concentrate our attention to the construction of the quantity
$G_3$ which is in front of the homogeneity of degree $8$ in $x$ and
$y$ in \eqref{eq49}. Writing the system \eqref{eq35}--\eqref{eq40}
in the matrix form
\begin{equation*}
\begin{gathered}
A_{3} B_{3}=C_{3},
\end{gathered}
\end{equation*}
we obtain
\begin{equation}
\begin{gathered}\label{eq55}
G_3=\frac{G_{3,i,j,k}+B_{3,i,j,k}b_i+D_{3,i,j,k}d_j+F_{3,i,j,k}f_k}{\sigma_{3,i,j,k}},\
(i=\overline{0,4};\ j=\overline{0,6};\ k=\overline{0,8}).
\end{gathered}
\end{equation}

Similarly to the previous case, we choose a comitant of the weight
$-1$ of the system $s(1,2)$ from \eqref{eq33} which contains as
semi-invariant the expression
$G_{3,2,j,k}+B_{3,2,j,k}b_2+D_{3,2,j,k}d_{j}+ +F_{3,2,j,k}f_k\
(k=\overline{0,8})$, and we find that it belongs to the space
\begin{equation*}
\begin{gathered}
S_{1,2}^{(8,37,6)}.
\end{gathered}
\end{equation*}

\section{General type of comitants which have as coefficients
expressions with generalized focal pseudo-quantities of the system
(34)}

\ \ \ \ Let's consider the extension of the system
\eqref{eq35}--\eqref{eq40} obtained from the identity \eqref{eq29}
for the system \eqref{eq33} and the function \eqref{eq34} which
contains the quantity $G_k$, which we write in a matrix form as
follows $A_kB_k=C_k$. We denote by $m_{G_{k}}$ the number of
equations and by $n_{G_{k}}$ the number of unknowns of this system.
Observe that these numbers can be written as
\begin{equation*}
\begin{gathered}
m_{G_{k}}=\underbrace{4+5}_{G_{1}}+\underbrace{6+7}_{G_{2}}+\underbrace{8+9}_{G_{3}}+\cdots+\underbrace{(2k+2)+(2k+3)}_{G_{k}},\
(k=1,2,3,...),\\
n_{G_{k}}=\underbrace{4+6}_{G_{1}}+\underbrace{6+8}_{G_{2}}+\underbrace{8+10}_{G_{3}}+\cdots+\underbrace{(2k+2)+(2k+4)}_{G_{k}}.
\end{gathered}
\end{equation*}
Hence we obtain
\begin{equation}
\begin{gathered}\label{eq56}
m_{G_{k}}=k(2k+7), n_{G_{k}}=m_{G_{k}}+k>m_{G_{k}}.
\end{gathered}
\end{equation}

Similarly to the previous cases, from this system we have
\begin{equation}
\begin{gathered}\label{eq57}
G_k=\frac{G_{k,i_{1},i_{2},\ldots, i_{k}}+B_{k,i_{1},i_{2},\ldots,
i_{k}}b_{i_1}+\cdots+Z_{k,i_{1},i_{2},\ldots,
i_{k}}z_{i_{k}}}{\sigma_{k,i_{1},i_{2},\ldots, i_{k}}},
\end{gathered}
\end{equation}

Now it is important to determine the degree of the polynomial
$G_{k,i_{1},i_{2},\ldots, i_{k}}$ in coefficients of the
differential system \eqref{eq33}.

Observe that the degree of non-zero polynomial coefficient of
$G_{i}\ (i=\overline{1,k)}$ in coefficients of the system
\eqref{eq33} in the matrix of Cramer's determinant of the order
$m_{G_{k}}$, when the column corresponding to the last quantity $G_{k}$ is replaced with the column corresponding to free
members, forms the following diagram (the last quantity $G_k$ has the degree $2$ according to the substitution):
\begin{equation*}
\begin{CD}
\ \ G_{1},G_{2},G_{3},...,G_{k-1},G_{k}.\\
\ \downarrow\ \ \ \downarrow\ \ \ \ \downarrow\ \ \ \ \ \ \ \  \downarrow\ \ \ \ \ \ \downarrow \\
\ 2\ \ \ 3\ \ \ \ 4\ \ \ \ \ \  \ \ k\ \ \ \ \ \ 2
\end{CD}
\end{equation*}
Then the degree of the polynomial  $G_{k,i_{1},i_{2},\ldots, i_{k}}$
in coefficients of the system \eqref{eq33}, denoted by $N_{G_{k}}$,
can be written as
\begin{equation*}
N_{G_{k}}=m_{G_{k}}-k+\frac{k(k+1)}{2}+1,
\end{equation*}
hence according to \eqref{eq56} we have
\begin{equation}
\begin{gathered}\label{eq58}
N_{G_{k}}=\frac{1}{2}(5k^2+13k+2).
\end{gathered}
\end{equation}

It is the degree of homogeneity of $G_{k,i_{1},i_{2},\ldots, i_{k}}$
in coefficients of the linear and the quadratic parts of the
differential system \eqref{eq33} which is contained in a polynomial
of the type $(d)=$ $=(\delta,d_1,d_2)$. Since $\delta=2(k+1)$ and
$d_2=2k$, then $d_1=N_{G_k}-2k$. So we obtain that a comitant of the
weight $-1$ of the system $s(1,2)$ from \eqref{eq33}, containing the
semi-invariant $G_{k,i_{1},i_{2},\ldots, i_{k}}+$
$+B_{k,i_{1},i_{2},\ldots,
i_{k}}b_{i_1}+\cdots+Z_{k,i_{1},i_{2},\ldots, i_{k}}z_{i_{k}}$,
which corresponds to the quantity $G_{k}$ for $k=1,2,3,...$, belongs
to the type
\begin{equation}\label{eq59}
\left( 2(k+1),\frac{1}{2}(5k^2+9k+2),2k \right),
\end{equation}
where $2(k+1)$ is the degree of homogeneity of the comitant in phase
variables $x,y$; $\displaystyle{\frac{1}{2}}(5k^2+9k+2)$ is the
degree of homogeneity of the comitant in coefficients of the linear
part $c,d,e,f$ and $2k$ is the degree of homogeneity of the comitant
in coefficients of the quadratic part of the system \eqref{eq33}.

Hereafter the expressions $G_{k,i_{1},i_{2},\ldots, i_{k}}$, which
determine the types of comitants \eqref{eq59} corresponding to the
quantity $G_{k}\ (k=1,2,3,...)$ will be called {\it the defining
focal quantities}. The comitants of the type \eqref{eq59} for
$k=1,2,3,...$ which contains as the coefficients expressions with the
generalized focal pseudo-quantities
\begin{equation*}
G_{k,i_{1},i_{2},\ldots, i_{k}}+B_{k,i_{1},i_{2},\ldots,
i_{k}}b_{i_1}+\cdots+Z_{k,i_{1},i_{2},\ldots, i_{k}}z_{i_{k}}.
\end{equation*}
will be called {\it the comitants associated to generalized focal
pseudo-quantities}.

For $G_{0}$ from \eqref{eq31}, which for the system $s(1,2)$ from
\eqref{eq33} has the type $(0,1,0)$, we retain the name {\it a null
focal pseudo-quantity}.

The space of comitants of the system $s(1,2)$ from \eqref{eq33},
corresponding to the type \eqref{eq59}, will be denoted by
\begin{equation}\label{eq60}
S_{1,2}^{(2(k+1),\frac{1}{2}(5k^2+9k+2),2k)}.
\end{equation}

\section{Comitants which have as coefficients expressions with
generalized focal pseudo-quantities of the system
$\boldsymbol{s}$(1,2,3)}

\ \ \ \ Let us consider the system $s(1,2,3)$ of the form
\begin{equation}\label{eq062}
\begin{gathered}
\dot{x}=cx+dy+gx^2+2hxy+kx^2+px^3+3qx^2y+3rxy^2+sy^3,\\
\dot{y}=ex+fy+lx^2+2mxy+ny^2+tx^3+3ux^2y+3cxy^2+wy^3
\end{gathered}
\end{equation}
with finitely determined Sibirsky graded algebra of unimodular
comitants $S_{1,2,3}$ [7], for which the function \eqref{eq30} will
be write in the form \eqref{eq34}, where $k_2$ is from \eqref{eq25}
and $a_0,a_1,...,f_7,f_8,...,G_1,G_2,...$ are unknowns. Similarly as
in the Sections 6 and 7 for determining the quantity $G_1$ we write
the equations in which splits the identity \eqref{eq29} in the case
of the system \eqref{eq062} in the matrix form
\begin{equation}\label{eq063}
\begin{gathered}
\widetilde{A}_1\widetilde{B}_1=\widetilde{C}_1,
\end{gathered}
\end{equation}
where
{\scriptsize
$$
\begin{gathered}
\widetilde{A}_1\!=\!\left(\!\!\!
\begin{array}{cccccccccc}
3 c& 3 e& 0& 0& 0& 0& 0& 0& 0& 0\\
 3 d& 6 c\! +\! 3 f& 6 e& 0& 0& 0& 0&
0& 0& 0\\
0& 6 d& 3 c\! + \!6 f& 3 e& 0& 0& 0& 0& 0& 0\\
0& 0& 3 d& 3 f& 0& 0& 0& 0& 0& 0\\
3 g& 3 l& 0& 0& 4 c& 4 e& 0& 0& 0& -e^2\\
6 h& 6 g\! +\! 6 m& 6 l& 0& 4 d& 12 c\! +\! 4 f& 12 e& 0& 0&
   2 c e\! -\! 2 e f \\
3 k& 12 h\! + \!3 n& 3 g\! +\! 12 m& 3 l& 0& 12 d& 12 c \!+ \!12 f& 12 e&
   0& -c^2 \!+ \!2 d e \!+ \!2 c f\! - \!f^2 \\
0& 6 k& 6 h\! +\! 6 n& 6 m& 0& 0& 12 d& 4 c\! +\! 12 f&
   4 e& -2 c d\! +\! 2 d f \\
0& 0& 3 k& 3 n& 0& 0& 0& 4 d& 4 f& -d^2\\
 \end{array}\!\!\!
\right),
\end{gathered}
$$
}
{\footnotesize
\begin{equation}\label{eq064}
\begin{gathered}
\widetilde{B}_1=\left(%
\begin{array}{c}
  a_0 \\
  a_1 \\
  a_2 \\
  a_3 \\
  b_0 \\
  b_1 \\
  b_2 \\
  b_3 \\
  b_4 \\
  G_1 \\
\end{array}%
\right),\quad \widetilde{C}_1=\left(%
\begin{array}{c}
2 e g - c l + f l\\
 -c g + f g + 4 e h - 2 d l - 2 c m + 2 f m\\
  -2
c h + 2 f h + 2 e k - 4 d m - c n + f n\\
 -c k + f k - 2 d n\\
  2 e p- c t + f t\\
 -c p + f p + 6 e q - 2 d t - 3 c u + 3 f u\\
  -3 c q + 3 f q + 6 e r - 6 d u - 3 c v + 3 f v\\
  -3 c r + 3 f r + 2 e s - 6 d v - c w + f w\\
 -c s + f s - 2 d w\\
\end{array}%
\right)
\end{gathered}
\end{equation}
}
For each fixed  $i\in\{0,1,...,4\}$ using the Cramer's rule from the
system \eqref{eq063} we find
\begin{equation}\label{eq065}
\begin{gathered}
\widetilde{G}_1=\frac{\widetilde{G}_{1,i}+\widetilde{B}_{1,i}b_i}{\widetilde{\sigma}_{1,i}},
\end{gathered}
\end{equation}
where
$\widetilde{G}_{1,i},\widetilde{B}_{1,i},\widetilde{\sigma}_{1,i}$
are polynomials in the coefficients of the system \eqref{eq062} and
$b_i$ are undetermined coefficients of the function $U(x,y)$ from
\eqref{eq34}.

By studying the matrices \eqref{eq063}--\eqref{eq064} of the system
\eqref{eq062} we conclude that the focal pseudo-quantity
$\widetilde{G}_{1,i}$ for fixed $i$ from \eqref{eq065} can be write as
\begin{equation}\label{eq066}
\begin{gathered}
\widetilde{G}_{1,i}=\widetilde{G}_{1,i}^{\prime}+\widetilde{G}_{1,i}^{\prime\prime},\
(i=0,1,2,3,4),
\end{gathered}
\end{equation}
where $\widetilde{G}_{1,i}^{\prime}$ (respectively
$\widetilde{G}_{1,i}^{\prime\prime})$ are homogeneous polynomials of
degree $8$ (respectively $9$) in coefficients of the linear part and
of degree $2$ in the coefficients of the quadratic part (respectively of
degree $1$ in the coefficients of the cubic part) of the differential
system \eqref{eq062}.

Using here the operators \eqref{eq17} of Lie algebra $L_4$ from [7]
for the system \eqref{eq062} we construct the corresponding operators which we denote respectively by
$\mathcal{X}_1,\mathcal{X}_2,\mathcal{X}_3,\mathcal{X}_4$.  Applying these operators under the
expressions from \eqref{eq066} we find
\begin{equation*}
\begin{gathered}
\mathcal{X}_1(\widetilde{f}_4^{\prime})=\mathcal{X}_4(\widetilde{f}_4^{\prime})=\widetilde{f}_4^{\prime},\
\mathcal{X}_2(\widetilde{f}_4^{\prime})=\mathcal{X}_3(\widetilde{f}_4^{\prime})=0,\\
\mathcal{X}_1(\widetilde{f}_4^{\prime\prime})=\mathcal{X}_4(\widetilde{f}_4^{\prime\prime})=\widetilde{f}_4^{\prime\prime},\
\mathcal{X}_2(\widetilde{f}_4^{\prime\prime})=\mathcal{X}_3(\widetilde{f}_4^{\prime\prime})=0,
\end{gathered}
\end{equation*}
where
\begin{equation}\label{eq067}
\begin{gathered}
\widetilde{f}_4^{\prime}(x,y)=\widetilde{G}_{1,0}^{\prime}x^4-4\widetilde{G}_{1,1}^{\prime}x^3y+2\widetilde{G}_{1,2}^{\prime}x^2y^2+4\widetilde{G}_{1,3}^{\prime}xy^3-\widetilde{G}_{1,4}^{\prime}y^4,\\
\widetilde{f}_4^{\prime\prime}(x,y)=\widetilde{G}_{1,0}^{\prime\prime}x^4-4\widetilde{G}_{1,1}^{\prime\prime}x^3y+2\widetilde{G}_{1,2}^{\prime\prime}x^2y^2+4\widetilde{G}_{1,3}^{\prime\prime}xy^3-\widetilde{G}_{1,4}^{\prime\prime}y^4,
\end{gathered}
\end{equation}
are comitants of the weight  $-1$ of the system \eqref{eq062} and
$\widetilde{G}_{1,i}^{\prime}$, $\widetilde{G}_{1,i}^{\prime\prime}$
are from \eqref{eq065}.

According to the above mentioned and \eqref{eq4} the given comitants
\eqref{eq067} belongs  respectively to the linear spaces
\begin{equation}\label{eq068}
\begin{gathered}
S_{1,2,3}^{(4,8,2,0)},\ S_{1,2,3}^{(4,9,0,1)},
\end{gathered}
\end{equation}
which are components of Sibirsky graded algebra of comitants
$S_{1,2,3}$ for the system \eqref{eq062}.

Taking into account \eqref{eq065} for $b_i=0\ (i=\overline{0,4})$ on
the variety $\mathcal{V}$ from \eqref{eq26} and also \eqref{eq066},
\eqref{eq067} we find out that  the first focal quantity $L_1$ of the
system \eqref{eq062} is related to the comitants \eqref{eq067} as
follows
\begin{equation*}
\begin{gathered}
\left[\widetilde{f}_4^{\prime}(x,y)+\widetilde{f}_4^{\prime\prime}(x,y)\right]\left|\right
._{\mathcal{V}}=8L_1(x^2+y^2)^2,
\end{gathered}
\end{equation*}
where
\begin{equation*}
\begin{gathered}
L_1=\frac{1}{4}\left\{[g(l-h)-k(h+n)+m(l+n)]-3[p+r+u+v]\right\}.
\end{gathered}
\end{equation*}

Similarly to the previous case, for determining the quantity $G_2$
for the system \eqref{eq062}, from the identity \eqref{eq29} we
obtain the following equation in the matrix form
\begin{equation}\label{eq069}
\begin{gathered}
\widetilde{A}_2\widetilde{B}_2=\widetilde{C}_2.
\end{gathered}
\end{equation}
For each fixed $i\in\{0,1,...,4\},\ j\in\{0,1,2,...,6\}$
we find the expression
\begin{equation}\label{eq070}
\begin{gathered}
\widetilde{G}_2=\frac{\widetilde{G}_{2,i,j}+\widetilde{B}_{2,i,j}b_i+\widetilde{D}_{2,i,j}d_j}{\widetilde{\sigma}_{2,i,j}}.
\end{gathered}
\end{equation}

By studying the matrix equality \eqref{eq069} we find that the focal
pseudo-quantity from \eqref{eq070} can be written in the form of
homogeneity of degree $24$ that can be represented in the form
\begin{equation}\label{eq071}
\begin{gathered}
\widetilde{G}_{2,i,j}=\widetilde{G}_{2,i,j}^{\prime}+\widetilde{G}_{2,i,j}^{\prime\prime}+\widetilde{G}_{2,i,j}^{\prime\prime\prime},
\end{gathered}
\end{equation}
where $ \widetilde{G}_{2,i,j}^{\prime},\
\widetilde{G}_{2,i,j}^{\prime\prime}$ and $
\widetilde{G}_{2,i,j}^{\prime\prime\prime},$ are homogeneity of
the type \eqref{eq4} respectively of the form $(0,20,4,0)$, $(0,21,2,1)$ and
$(0,22,0,2).$ We note that on the variety $\mathcal{V}$ from
\eqref{eq26} for the system \eqref{eq062} the quantities
$\widetilde{G}_{2,2,j}\ (j=\overline{0,6})$ have the expressions
\begin{equation*}
\begin{gathered}
\widetilde{G}_{2,2,j}|_{\mathcal{V}}=2304L_2,\  (j=0,2,4,6),\
\widetilde{G}_{2,2,j}|_{\mathcal{V}}=0,\ (j=1,3,5).
\end{gathered}
\end{equation*}

On the other hand, the second focal quantity $L_2$ of the system
\eqref{eq062} can be written with the terms from \eqref{eq071} as
follows
\begin{equation*}
\begin{gathered}
24L_2=\widetilde{G}_{2,2,j}^{\prime}|_{\mathcal{V}}+\widetilde{G}_{2,2,j}^{\prime\prime}|_{\mathcal{V}}+\widetilde{G}_{2,2,j}^{\prime\prime\prime}|_{\mathcal{V}},\
(j=0,2,4,6),
\end{gathered}
\end{equation*}
where
\begin{equation*}
\begin{gathered}
\widetilde{G}_{2,2,j}^{\prime}|_{\mathcal{V}}=4(62 g^3 h - 2 g h^3 + 95 g^2 h k - 2 h^3 k + 38 g h k^2 + 5 h k^3 -
  62 g^3 l + 27 g h^2 l -\\- 39 g^2 k l + 29 h^2 k l - 15 g k^2 l -
  8 g h l^2 + 15 h k l^2 - 5 g l^3 + 53 g^2 h m + 66 g h k m +\\+
  13 h k^2 m - 127 g^2 l m - 6 h^2 l m - 68 g k l m - 15 k^2 l m -
  13 h l^2 m - 5 l^3 m +\\+ 6 g h m^2 + 6 h k m^2 - 63 g l m^2 -
  29 k l m^2 + 2 l m^3 + 6 g^3 n + 61 g h^2 n + 72 g^2 k n +\\+
  63 h^2 k n + 33 g k^2 n + 5 k^3 n - 10 g h l n + 68 h k l n -
  33 g l^2 n + 15 k l^2 n - 72 g^2 m n -\\- 6 h^2 m n + 10 g k m n +
  8 k^2 m n - 66 h l m n - 38 l^2 m n - 61 g m^2 n - 27 k m^2 n +\\+
  2 m^3 n + 72 g h n^2 + 127 h k n^2 - 72 g l n^2 + 39 k l n^2 -
  53 h m n^2 - 95 l m n^2 - \\-6 g n^3 + 62 k n^3 - 62 m n^3),
\end{gathered}
\end{equation*}
\begin{equation*}
\begin{gathered}
\widetilde{G}_{2,2,j}^{\prime\prime}|_{\mathcal{V}}=-2(186 g^2 p + 10 h^2 p + 117 g k p + 45 k^2 p + 59 h l p + 15 l^2 p +
 159 g m p + \\+75 k m p + 18 m^2 p + 143 h n p + 89 l n p + 196 n^2 p -
 69 g h q - 57 h k q + 69 g l q +\\+ 12 k l q + 9 l m q + 60 g n q +
 3 k n q + 21 m n q + 168 g^2 r - 6 h^2 r + 69 g k r + 15 k^2 r +\\+
 87 h l r + 45 l^2 r + 123 g m r + 39 k m r + 18 m^2 r + 171 h n r +
 129 l n r + 222 n^2 r -\\- 13 g h s - 17 h k s - 15 g l s - 16 h m s -
 15 l m s - 16 g n s - 17 k n s - 19 m n s -\\- 19 g h t - 15 h k t -
 17 g l t - 16 h m t - 17 l m t - 16 g n t - 15 k n t - 13 m n t +\\+
 222 g^2 u + 18 h^2 u + 129 g k u + 45 k^2 u + 39 h l u + 15 l^2 u +
 171 g m u + 87 k m u -\\- 6 m^2 u + 123 h n u + 69 l n u + 168 n^2 u +
 21 g h v + 9 h k v + 3 g l v + 12 k l v - 57 l m v +\\+ 60 g n v +
 69 k n v - 69 m n v + 196 g^2 w + 18 h^2 w + 89 g k w + 15 k^2 w +
 75 h l w +\\+ 45 l^2 w +143 g m w + 59 k m w + 10 m^2 w + 159 h n w +
 117 l n w + 186 n^2 w),
\end{gathered}
\end{equation*}
\begin{equation*}
\begin{gathered}
\widetilde{G}_{2,2,j}^{\prime\prime\prime}|_{\mathcal{V}}=-9(11 p q + 15 q r - 5 p s - r s + p t + 5 r t + 3 q u - 5 s u + t u -
 7 p v - 3 r v - \\-15 u v + 7 q w - s w + 5 t w - 11 v w).
\end{gathered}
\end{equation*}

Similarly to the technique described in the Sections 6 and 7 we choose a
comitant of the weight $-1$ of the system $s(1,2,3)$ from \eqref{eq062}
which contains as a semi-invariant the expression
$\widetilde{G}_{2,i,j}+\widetilde{B}_{2,i,j}b_i+\widetilde{D}_{2,i,j}d_j$.
According to the decomposition \eqref{eq071} and the types shown
below we find that this comitant is a sum of comitants belonging to
the spaces
\begin{equation}\label{eq072}
S_{1,2,3}^{(6,20,4,0)},\ S_{1,2,3}^{(6,21,2,1)},\
S_{1,2,3}^{(6,22,0,2)}.
\end{equation}

Following this process with the help of the matrix equation
\begin{equation*}
\begin{gathered}
\widetilde{A}_3\widetilde{B}_3=\widetilde{C}_3
\end{gathered}
\end{equation*}
for each fixed $i\in\{0,1,...,4\},\ j\in\{0,1,...,6\},\
k\in\{0,1,...,8\}$ we obtain
\begin{equation}\label{eq073}
\begin{gathered}
\widetilde{G}_3=\frac{\widetilde{G}_{3,i,j,k}+\widetilde{B}_{3,i,j,k}b_i+\widetilde{D}_{3,i,j,k}d_j+\widetilde{F}_{3,i,j,k}f_j}{\widetilde{\sigma}_{3,i,j,k}}.
\end{gathered}
\end{equation}

Similarly to the previous case we find that the focal
pseudo-quantity $\widetilde{G}_{3,i,j,k}$ splits into a sum of
four terms of the same degree $43$ in the coefficients of the system
\eqref{eq062}, which according to \eqref{eq4} belongs to the types
$(0,37,6,0)$, $(0,38,4,1)$, $(0,39,2,2)$ and $(0,40,0,3)$. Then it
results that the comitant of the weight $-1$ having as a semi-invariant
one of the expressions \eqref{eq073} consists of the sum of
comitants of the system \eqref{eq062} which belongs to the spaces
\begin{equation}\label{eq074}
S_{1,2,3}^{(8,37,6,0)},\ S_{1,2,3}^{(8,38,4,1)},\
S_{1,2,3}^{(8,39,2,2)},\ S_{1,2,3}^{(8,40,0,3)}.
\end{equation}

Following this process we obtain the sequence of linear
spaces \eqref{eq068}, \eqref{eq072}, \eqref{eq074} etc. of comitants of the
system \eqref{eq062}. It remain to underline that the generalized focal pseudo-quantities corresponding to $G_k$ of the given system is exactly a sum of coefficients of these comitants.

It is not difficult to deduce that the generic formula of the types
of the comitants in which the generalized focal pseudo-quantities
corresponding to $G_k,\ (k=1,2,3,...)$ splits as a sum, has the
form:
\begin{equation*}
\begin{gathered}
\left(2(k+1),\frac{1}{2}(5k^2+9k+2)+i,2(k-i),i\right),\ (i=\overline{0,k}).
\end{gathered}
\end{equation*}

\section{Graded algebra of comitants whose spaces contain
comitants associated to generalized focal pseudo-quan-tities of the
system (34) and (62)}

\ \ \ \ Thus we obtain for the system \eqref{eq33} the set of spaces
of center-affine (unimodular) comitants
\begin{equation}\label{eq61}
\mathbb{R}=S_{1,2}^{(0,0,0)},\ S_{1,2}^{(0,1,0)},\
S_{1,2}^{(4,8,2)},\
S_{1,2}^{(6,20,4)},...,S_{1,2}^{(2(k+1),\frac{1}{2}(5k^2+9k+2),2k)},...\subset
S_{1,2},
\end{equation}
were $S_{1,2}$ is Sibirsky graded algebra of the system
\eqref{eq33}.

Let's consider the graded algebra $S_{1,2}^{\prime}$, generated by
the space $S_{1,2}^{(\delta^\prime,d_{1}^\prime,d_2^\prime)}$ from
\eqref{eq61}, which can be written as
\begin{equation}
\begin{gathered}\label{eq62}
S_{1,2}^{\prime}=\bigoplus_{(d^{\prime})}S_{1,2}^{(d^{\prime})}.
\end{gathered}
\end{equation}
Here $S_{1,2}^{(d^{\prime})}$ denote linear spaces, contained in
$S_{1,2}^{(\delta^\prime,d_{1}^\prime,d_2^\prime)}$ for all
$(d^{\prime})$, as well as the spaces from $S_{1,2}$ which contains all possible products of spaces \eqref{eq61}.

Since the algebra $S_{1,2}^{\prime}$ is a graded subalgebra of the
algebra $S_{1,2}$ for the system \eqref{eq33}, according to
Proposition 4, we obtain that for the Krull dimensions of these
algebras the following inequality takes place:
\begin{equation}
\begin{gathered}\label{eq63}
\varrho(S_{1,2}^{\prime})\leq\varrho(S_{1,2}).
\end{gathered}
\end{equation}

Taking into account this inequality, and the fact that from [7] we
have $\varrho(S_{1,2})=9$, according to Definition 2, we have

\noindent\textbf{Lemma 1}. {\it The maximal number of algebraically
independent generalized focal pseudo-quantities in the center-focus
 problem   for    the  system \eqref{eq33} does  not  exceed $9$}.

According to the equalities \eqref{eq31}, \eqref{eq47}, \eqref{eq54}
etc. and the conclusion, resulting from Proposition 5, that the
number of algebraically independent focal quantities $L_k\
(k=\overline{0,\infty})$ can not exceed the maximal number of
algebraically independent generalized focal pseudo-quantities, using
Lemma 1, we have

\noindent\textbf{Theorem 2}. {\it The maximal number of
algebraically independent focal quantities of the system
\eqref{eq33} on the variety $\mathcal{V}$ from \eqref{eq26} or,
equivalently, from \eqref{eq28}, that take part in solving the
center-focus problem, does not exceed $9$}.

With the help of Hilbert series of the algebras $S_{1,2}$,
$S_{1,2}^{\prime}$, $SI_{1,2}$ [23] and Remark 3 it can be shown
that the predicted upper bound of algebraically independent focal
quantities of the system \eqref{eq33} on the variety $\mathcal{V}$
from \eqref{eq26} (\eqref{eq28}) can be much less than $9$, and can
be equal to $7$ or, may be, even $5$.

We note that the similar studies that for the system $s(1,2)$ from
\eqref{eq33} were realized in the works [27,28,31] for the systems
$s(1,3),\ s(1,4),\ s(1,5)$ respectively. This scenario is confirmed
for the system $s(1,2,3)$ by studies in the case 9 which allow to
form the algebra $S_{1,2,3}^{\prime}$ with the same properties as
the algebra $S_{1,2}^{\prime}$.

\section{Main results}

\ \ \ \ Similarly to the considered cases, for any system
$s(1,m_1,m_2,...,m_{\ell})$ from \eqref{eq1} we have that the
algebras similar to the obtained in the above mentioned examples
satisfy the inclusion condition
\begin{equation*}
\begin{gathered}
S_{1,m_1,m_2,...,m_{\ell}}^{\prime}\subset
S_{1,m_1,m_2,...,m_{\ell}},
\end{gathered}
\end{equation*}
hence according to Proposition 5, for their Krull dimensions we have
\begin{equation}
\begin{gathered}\label{eq64}
\varrho(S_{1,m_1,m_2,...,m_{\ell}}^{\prime})\leq\varrho(S_{1,m_1,m_2,...,m_{\ell}}).
\end{gathered}
\end{equation}
By the formula \eqref{eq18} we obtain
\begin{equation}
\begin{gathered}\label{eq65}
\varrho(S_{1,m_1,m_2,...,m_{\ell}})=2\left(\sum_{i=1}^{\ell}m_i+\ell\right)+3.
\end{gathered}
\end{equation}
Analogously to Lemma 1 and other considered examples, with the help
of \eqref{eq64} and \eqref{eq65} it can be shown that the following
lemma is true:

\noindent\textbf{Lemma 2}. {\it The maximal number of algebraically
independent generalized focal pseudo-quantities in the center-focus
problem for the system \eqref{eq1} does not exceed the number from
\eqref{eq65}}.

\noindent\textit{Remark} 8.  According to the Remarks 5 and 6 and
formulation of center-focus problem given in Section 5, as well as
the identities \eqref{eq31} we can say that the generalized focal
pseudo-quantities, being semi-invariants in above mentioned
comitants, have as projections on the variety  $\mathcal{V}$ from
\eqref{eq26} (\eqref{eq28}) the focal quantities $L_k\
(k=1,2,3,...).$

From identity \eqref{eq29} and Lyapunov's function \eqref{eq34} it
results that for any system $s(1,m_1,$ $m_2,...,m_{\ell})$ we can
write the identities of the type \eqref{eq57} for quantities $G_k$
$(k=1,2,3,...)$, which have as numerators the generalized focal
pseudo-quantities. Using these quantities and the operator $D_3$
from \eqref{eq17} of system $s(1,m_1,m_2,...,m_{\ell})$ we can
determine comitants of the given system, having as coefficients the
above mentioned focal pseudo-quantities.

According to the Remark 8 we conclude that the following statement
take place:

\noindent\textbf{Theorem 3}. {\it The maximal number of
algebraically independent focal quantities of the system \eqref{eq1}
on the variety $\mathcal{V}$ from \eqref{eq26} or, equivalently,
from \eqref{eq28}, that take part in solving the center-focus
problem does not exceed the number from \eqref{eq65}}.

We recall that in the introduction it was told that for the systems
$s(1,2)$ and $s(1,3)$ the number of essential conditions for center
$\omega=3$ and $5$, respectively, but for the system $s(1,2,3)$
there is an assumption that $\omega\leq13$.

From Theorem 3 we obtain that the maximal number of algebraically
independent focal quantities for the system $s(1,2)$ does not exceed
$9$, for $s(1,3)$ does not exceed $11$, and for $s(1,2,3)$ does not
exceed $17$.

These arguments and Proposition 5 with $\mathcal{V}$ from
\eqref{eq26} or, equivalently, from \eqref{eq28}, and the defined
above algebra $S_{1,m_1,m_2,...,m_{\ell}}^{\prime}$ suggest that is
true

\noindent\textbf{The main hypothesis}. {\it The number $\omega$ of
essential conditions for center  from \eqref{eq3} which solve the
center-focus problem for the system \eqref{eq1}, having at the
origin of coordinates a singular point of the second type, does not
exceed the number from \eqref{eq65}}.

\noindent\textit{Remark} 9. The equality \eqref{eq65} shows that the
quantity $\varrho$ is equal to the number of coefficients of the
right parts of the system \eqref{eq1} minus one.

Besides [23], the authors have published their vision of the
center-focus problem in the theses [24-33].

\end{document}